\documentclass{article}
\usepackage{amsmath}
\usepackage{amsfonts}
\usepackage{amssymb}	
\usepackage{color}
\newtheorem{theo}{Theorem}[section]
\newtheorem{lem}[theo]{Lemma}
\newtheorem{cor}[theo]{Corollary}
\newtheorem{rem}[theo]{Remark}

\newcommand{\mysection}[1]{\section{#1} \setcounter{equation}{0}}
\newcommand{\proof}{{\sc Proof.} \quad}

\newcommand{\R}{\mathbb{R}}
\newcommand{\N}{\mathbb{N}}
\newcommand{\be}{\begin{equation} \label}
\newcommand{\ee}{\end{equation}}
\newcommand{\bes}{\begin{equation} \begin{array}{c} \label}
\newcommand{\ees}{\end{array} \end{equation}}
\newcommand{\bea}{\begin{eqnarray}\label}
\newcommand{\eea}{\end{eqnarray}}
\newcommand{\beas}{\begin{eqnarray} \begin{array}{rcl} \label}
\newcommand{\eeas}{\end{array} \end{eqnarray}}
\newcommand{\bas}{\begin{eqnarray*}}\newcommand{\eas}{\end{eqnarray*}}
\newcommand{\bass}{\begin{eqnarray*} \begin{array}{rcl}}
\newcommand{\eass}{\end{array} \end{eqnarray*}}
\newcommand{\basss}{\begin{eqnarray*} \begin{array}{c}}
\newcommand{\easss}{\end{array} \end{eqnarray*}}
\newcommand{\qed}{{}\hfill $\square$ \\}
\newcommand{\bit}{\begin{itemize}}
\newcommand{\eit}{\end{itemize}}
\newcommand{\nn}{\nonumber}
\newcommand{\eps}{\varepsilon}
\newcommand{\abs}{\\[3mm]}

\newcommand{\pO}{\partial\Omega}

\newcommand{\io}{\int_\Omega}

\newcommand{\F}{{\cal F}}
\newcommand{\D}{{\cal D}}
\newcommand{\B}{{\cal B}}
\newcommand{\set}{{\cal S}(m,M,B,\kappa)}
\newcommand{\tm}{T_{max}(u_0,v_0)}

\begin{document}
\title{Finite-time blowup and global-in-time unbounded solutions to a parabolic-parabolic quasilinear Keller-Segel system in higher dimensions}
\author{
Tomasz Cie\'{s}lak \\
{\small Institut f\"ur Mathematik, Universit\"at Z\"{u}rich, Winterthurerstrasse 190, 8057 Z\"{u}rich, Switzerland} \\
{\small Instytut Matematyczny PAN, \'Sniadeckich 8, 00-956 Warsaw, Poland}\\
{\small E-Mail: T.Cieslak@impan.pl}
\and 
Christian  Stinner\\
{\small Institut f\"ur Mathematik, Universit\"at Z\"{u}rich, Winterthurerstrasse 190, 8057 Z\"{u}rich, Switzerland} \\
{\small E-Mail: christian.stinner@math.uzh.ch}
}
\date{}
\maketitle
\begin{abstract}
  \noindent
  In this paper we consider quasilinear Keller-Segel type systems of two kinds in higher dimensions. In the case of a nonlinear diffusion system we prove an optimal (with respect to possible nonlinear diffusions generating explosion in finite time of solutions) finite-time blowup result. In the case of a cross-diffusion system we give results which are optimal provided one assumes some proper non-decay of a nonlinear chemical sensitivity. Moreover, we show that once we do not assume  the above mentioned non-decay, our result cannot be as strong as in the case of nonlinear diffusion without nonlinear cross-diffusion terms. To this end we provide an example, interesting by itself, of global-in-time unbounded solutions to the nonlinear cross-diffusion Keller-Segel system with chemical sensitivity decaying fast enough, in a range of parameters in which there is a finite-time blowup result in a corresponding case without nonlinear cross-diffusion.  

\end{abstract}
\mysection{Introduction}\label{section1}
This work deals with radially symmetric nonnegative solution couples $(u,v)$ of the parabolic-parabolic 
Keller-Segel system
\be{0}
	\left\{ \begin{array}{ll}
	u_t= \nabla \cdot (\phi(u) \nabla u) - \nabla \cdot ( \psi(u) \nabla v ), & \; x\in\Omega, \ t>0, \\[2mm]
	v_t=\Delta v-v+u, & \; x\in\Omega, \ t>0, \\[2mm]
	\frac{\partial u}{\partial\nu}=\frac{\partial v}{\partial\nu}=0, & \; x\in\partial\Omega, \ t>0, \\[2mm]
	u(x,0)=u_0(x), \quad v(x,0)=v_0(x), & \; x\in\Omega,
	\end{array} \right.
\ee
in a ball $\Omega = B_R \subset \R^n$, where $n \ge 3$, $R>0$,
and the
initial data are supposed to satisfy $u_0 \in C^0(\bar\Omega)$ and $v_0\in W^{1,\infty}(\Omega)$ such that
$u_0 > 0$ and $v_0 >0$ in $\bar{\Omega}$.

Moreover, we assume that $\phi, \psi \in C^2([0,\infty))$ and that there is $\beta \in C^2 ([0,\infty))$ 
such that 
\begin{equation}\label{0.1}
 \phi(s) >0, \qquad \psi(s) = s \beta (s), \quad\mbox{and }\quad \beta (s) >0 \quad\mbox{for } s \in [0,\infty) 
\end{equation}
are satisfied.

Systems of this kind were introduced in \cite{KS:1} to describe the motion of cells which are diffusing and moving towards the gradient of a substance called chemoattractant, the latter being produced by the cells themselves. In particular, 
the essentiality of both nonlinear diffusion as well as nonlinear chemosensitivity were emphasized in \cite{hil:volume}
where it was explained that they can be used to model the so-called volume-filling effect. The Keller-Segel system has been studied extensively by many authors and the main issue of the investigation was chemotactic collapse of cells interpreted as finite-time blowup of the component $u$ of a solution to \eqref{0}. It is however worth to be underlined that despite the fact that the original Keller-Segel model was a system of parabolic equations the main results concerning the finite-time blowup of solutions to \eqref{0} were usually proved for its parabolic-elliptic simplification. There were a few methods introduced to investigate the phenomenon of finite-time explosion of solutions in that case. Two main methods among them being the change of variables leading to a reduction of the parabolic-elliptic simplification of \eqref{0} to a single equation obeying a maximum principle introduced in \cite{jl:expl} and the so-called moment method making strong use of the fact that the second equation of the parabolic-elliptic simplification of \eqref{0} is a Poisson equation, see  \cite{B},\cite{Na95}. Those two methods and their ramifications led to a variety of results concerning appearance of chemotactic collapse in both semilinear (i.e. $\phi=\psi\equiv 1$) and quasilinear Keller-Segel systems. In particular, there have been characterized values of initial mass distinguishing between finite-time blowup and global existence of bounded solutions to the two-dimensional semilinear version of \eqref{0} in both radially symmetric ond non-radial settings (see \cite{Na95}, \cite{Nagai}, \cite{NSY}). Moreover, it has been shown that in higher dimensions a finite-time blowup of solutions to the semilinear version of \eqref{0} can occur independently of the initial mass provided that the initial data are concentrated enough \cite{Na95}.  Finally, in the case of a quasilinear system, for any space dimension $n$ there have been identified critical nonlinearities such that if $\phi$ and $\psi$ satisfy the supercritical relation, then solutions to \eqref{0} stay bounded for any time while for those satisfying the subcritical relation solutions blow up in finite-time independently of the magnitude of initial mass provided the data are concentrated enough, see \cite{djiewin}.

However, all those results are available only for a parabolic-elliptic simplification of \eqref{0}. In the case of the original fully parabolic version the investigation of chemotactic collapse turned out to be a much more challenging issue. So far the only two existing results in the literature stating the occurrence of finite-time blowup of solutions to \eqref{0} are those in \cite{HV:1}, where an example of a special solution to the semilinear version of \eqref{0} in dimension $n=2$ blowing up in a finite-time is shown, and the result in \cite{clKS} where the explosion of solutions to the one-dimensional Keller-Segel system with appropriately weak diffusion of cells and sufficiently fast diffusion of chemoattractant is shown. The breakthrough has been made recently in \cite{win_bu}. Introducing a new method M. Winkler
shows there that in dimensions $n\geq 3$ generic solutions to the semilinear version of \eqref{0} blow up in finite time independently of the size of initial mass. In the present paper we generalize his method to the quasilinear case.
This way, to the best of our knowledge, we obtain a first result concerning a finite-time blowup of solutions to the fully parabolic quasilinear Keller-Segel system in higher dimensions. So far the only result in that direction was achieved in dimension $n=1$ and only for large initial masses in \cite{clKS}. Moreover, the result concerning a chemotactic collapse in the case where $\beta(u)\equiv 1$ is optimal. Namely, we show that in dimension $n\geq 3$, for $\phi(u)\leq Cu^{p}, p<1-\frac{2}{n}$ and some constant $C>0$, independently of the size of initial mass, one can find generic radially symmetric initial data leading to finite-time blowup (see Corollary \ref{cor1.35}). This result is optimal in view of the result in \cite{ss} guaranteeing global existence of bounded solutions to \eqref{0} with $\psi\equiv 1$ for $\phi$ satisfying $\phi(u)\geq Cu^{q}, q>1-\frac{2}{n}$ for some constant $C>0$. Moreover, in Corollary \ref{cor1.4} we prove that in the case of full nonlinear cross-diffusion we obtain a result at least as good as in the parabolic-elliptic case, compare \cite{djiewin}. Furthermore, it is also an optimal result for a fully parabolic problem when restricting ourselves to polynomial nonlinearities, see \cite{taowin_jde}. Theorem \ref{theo1}, which is our main achievement, shows that restricting ourselves to the case of $\psi(u)$ not decaying when $u$ is large, we obtain the result which is a counterpart of the existence of global-in-time solutions in \cite{ja2}. Finally, we show (see Theorem \ref{theo1.5}) that without assuming a lack of decay of $\psi(u)$ one cannot expect the existence of critical exponents distinguishing between boundedness of solutions and finite-time blowup. It turns out that the possible asymptotic behavior of solutions to the nonlinear cross-diffusion system \eqref{0} can be more complicated. We show that under the proper choice of parameters (corresponding to the choice of parameters which yield finite-time blowup in a semilinear case) one can construct global-in-time radially symmetric solutions admitting infinite-time blowup. This result seems to be quite interesting by itself since the phenomenon of infinite-time blowup does not seem to be that often met in parabolic equations.

To be more precise when formulating our finite-time blowup results we have to introduce the following notation.         Suppose that there exist positive constants $s_0$, $a$, and $b$ such that the functions
\begin{equation}\label{GH}
 G(s) := \int\limits_{s_0}^s \int\limits_{s_0}^\sigma \frac{\phi(\tau)}{\psi(\tau)} \; d\tau \, d \sigma, \quad s>0
 \qquad\mbox{and}\qquad 
 H(s) := \int\limits_0^s \frac{\sigma \phi(\sigma)}{\psi(\sigma)} \; d\sigma, \quad s \ge 0,
\end{equation}
fulfill
\begin{equation}\label{G1}
  G(s) \le a \, s^{2-\alpha}, \; s \ge s_0,  \qquad\mbox{with some } \alpha > \frac{2}{n},
\end{equation}
as well as
\begin{equation}\label{H1}
  H(s) \le \gamma \cdot G(s) + b(s+1), \; s > 0, \qquad\mbox{with some } \gamma \in \left( 0, \frac{n-2}{n} 
  \right).
\end{equation}
We remark that $H$ in \eqref{GH} is well-defined due to the positivity of $\beta$ in $[0,\infty)$.

It is well-known that the function 
\begin{equation}\label{F}
  \F(u,v):=\frac{1}{2} \io |\nabla v|^2 + \frac{1}{2} \io v^2 - \io uv + \io G(u)
\end{equation}
is a Liapunov functional for \eqref{0} with dissipation rate
\begin{equation}\label{D}
  \D(u,v):=\io v_t^2 + \io \psi(u) \cdot \Big| \frac{\phi (u)}{\psi(u)} \nabla u - \nabla v \Big|^2.
\end{equation}
More precisely, any classical solution to \eqref{0} satisfies
\begin{equation}\label{liapunov}
	\frac{d}{dt} \F(u(\cdot,t),v(\cdot,t)) = - \D (u(\cdot,t),v(\cdot,t)) 
	\qquad \mbox{for all } t \in (0,\tm),
\end{equation}
where $\tm \in (0,\infty]$ denotes the maximal existence time of $(u,v)$ (see \cite[Lemma~2.1]{win_mmas}). \abs

In order to prove our result of finite-time blowup, we need to impose the additional condition that there exists $c_0 >0$
such that
\begin{equation}\label{psi}
  \psi(s) \ge c_0 \, s, \quad s \ge 0,
\end{equation}
which in view of \eqref{0.1} means that $\beta (s) \ge c_0>0$ for $s \ge 0$.

Then we have the following result for blowup in finite time.
\begin{theo}\label{theo1}
  Suppose that $\Omega=B_R \subset \R^n$ with some $n\ge 3$ and $R>0$, assume that \eqref{G1}, \eqref{H1}, and \eqref{psi}
  are satisfied, and let $m>0$ and $A>0$ be given. 
  Then there exist positive constants $T(m,A)$ and $K(m)$ such that for any
  \begin{eqnarray}\label{t1.1}
	(u_0,v_0) \in \B (m,A) &:=& \bigg\{
	(u_0,v_0) \in C^0(\bar\Omega) \times W^{1,\infty}(\Omega) \ \bigg| \
	\mbox{$u_0$ and $v_0$ are radially symmetric} \nn \\ 
    & & \hspace*{5mm} \mbox{and positive in $\bar\Omega$, 
	$\io u_0=m$, $\|v_0\|_{W^{1,2}(\Omega)} \le A$,}\nn \\ 
    & & \hspace*{5mm} \mbox{and $\F(u_0,v_0) \le -K(m) \cdot (1+A^2)$} \bigg\},
  \end{eqnarray}
  the corresponding solution $(u,v)$ of \eqref{0} blows up at the finite time $\tm \in (0,\infty)$, where
  $\tm \le T(m,A)$.
\end{theo}

Moreover, the set $\B (m,A)$ has the following properties.
\begin{theo}\label{theo1.2}
  Let $\Omega=B_R \subset \R^n$ with some $n\ge 3$ and $R>0$, let $\B(m,A)$ be as defined in \eqref{t1.1}, and assume that 
  \eqref{G1} is fulfilled. 
  \begin{enumerate}
  \item[(i)] Then for any $m>0$ there exists $A>0$ such that $\B (m,A) \neq \emptyset$.
  \item[(ii)] Suppose that \eqref{G1} holds with some $\alpha > \frac{4}{n+2}$ and, moreover, let    
  $p \in (1,\frac{2n}{n+2})$ such that $p > 2-\alpha$.
  Then for any $m>0$ and $A>0$, the set $\B(m,A)$ 
  is dense in the space of all radially symmetric positive functions in $C^0(\bar\Omega) \times W^{1,\infty}(\Omega)$
  with respect to the topology in $L^p(\Omega) \times W^{1,2}(\Omega)$.
  In particular, given positive radial functions $(u_0,v_0) \in C^0(\bar\Omega) \times W^{1,\infty}(\Omega)$ and 
  $\eps>0$, there exist positive radial $(u_{0\eps},v_{0\eps}) \in C^0(\bar\Omega) \times W^{1,\infty}(\Omega)$
  such that 
  \bas
	\|u_{0\eps} - u_0\|_{L^p(\Omega)} + \|v_{0\eps}-v_0\|_{W^{1,2}(\Omega)} < \eps
  \eas
  and the solution $(u_\eps,v_\eps)$ of \eqref{0} with initial data $(u_{0\eps},v_{0\eps})$
  blows up in finite time.
  \end{enumerate}
\end{theo}
Furthermore, we state three Corollaries which cover interesting special cases. Corollary \ref{cor1.3} is an immediate consequence of Theorem~\ref{theo1} while 
Corollary \ref{cor1.35} follows since \eqref{H1} is satisfied which, in the case that $\phi$ is decreasing, is deduced in view of the possibility of choosing $s_0>e^{\frac{1}{\gamma}}$ and integration by parts, and, in case of $s^q / \phi(s) \to c$ as $s \to \infty$, is implied by \cite[Corollary~5.2]{win_mmas}. Moreover, Corollary~\ref{cor1.4} follows from Theorem~\ref{theo1}, because 
\cite[Corollary~5.2]{win_mmas} shows that the functions $\phi$ and $\psi$ given in Corollary~\ref{cor1.4} satisfy 
\eqref{G1} and \eqref{H1}. Corollary \ref{cor1.35} is optimal in view of the results given in \cite{ss}.
\begin{cor}\label{cor1.3}
  Assume that $\psi(s) =s$ for $s \ge 0$ and that \eqref{G1} and \eqref{H1} are fulfilled. Moreover, let 
  $\Omega=B_R \subset \R^n$ with some $n\ge 3$ and $R>0$,  and let $m>0$ and $A>0$ be given. 
  Then there exist positive constants $T(m,A)$ and $K(m)$ such that for any
  $(u_0,v_0) \in \B (m,A)$
  the corresponding solution $(u,v)$ of \eqref{0} blows up at the finite time $\tm \le T(m,A)$.
\end{cor}   
\begin{cor}\label{cor1.35}
  Assume that $\psi(s) =s$ for $s \ge 0$ and that 
  $\phi(s)\leq Cs^{q}$, $s \ge 1$, for some $q<1-\frac{2}{n}$ and $C>0$. Furthermore, suppose that either 
  $\phi$ is a decreasing function or that there exists $c>0$ such that $s^q / \phi(s) \to c$ as $s \to \infty$. Let
  $\Omega=B_R \subset \R^n$ with some $n\ge 3$ and $R>0$,  and let $m>0$ and $A>0$ be given. 
  Then there exist positive constants $T(m,A)$ and $K(m)$ such that for any
  $(u_0,v_0) \in \B (m,A)$
  the corresponding solution $(u,v)$ of \eqref{0} blows up at the finite time $\tm \le T(m,A)$.
\end{cor}
\begin{cor}\label{cor1.4}
  Assume that $\phi(s) = (s+1)^{-p}$ and $\psi(s) =s (s+1)^{q-1}$, $s \ge 0$, with $q \ge 1$ and $p \in \mathbb{R}$ such that 
  $p+q > \frac{2}{n}$. Moreover, let 
  $\Omega=B_R \subset \R^n$ with some $n\ge 3$ and $R>0$,  and let $m>0$ and $A>0$ be given. 
  Then there exist positive constants $T(m,A)$ and $K(m)$ such that for any
  $(u_0,v_0) \in \B (m,A)$
  the corresponding solution $(u,v)$ of \eqref{0} blows up at the finite time $\tm \le T(m,A)$.
\end{cor}   
In view of \cite{taowin_jde} the latter result is optimal in the case $q \ge 1$, while in view of \cite{win_mmas} it remains an interesting question whether Corollary~\ref{cor1.4} can be extended to the case $q<1$. In the following theorem, in particular, we provide a negative answer to this question. However, it still remains open to find critical exponents (if possible) 
distinguishing between finite- and infinite-time blowup of solutions when $q<1$. 

\begin{theo}\label{theo1.5}
Let $\Omega=B_R \subset \R^n$ with some $n\ge 2$ and $R>0$. Moreover, assume that $\lim_{s \to \infty} \phi (s) =0$, that 
there exists a positive constant $D>0$ such that for any $s>0$
\begin{equation}\label{balance}
\frac{\beta(s)}{\phi(s)}\leq D
\end{equation}
and that there exist constants $D_1>0$ and $\gamma_1>n$ such that for any $s>0$
\begin{equation}\label{decay}
\beta(s)\leq D_1s^{-\gamma_1}.
\end{equation} 
Assume also that \eqref{G1} and \eqref{H1} hold. Then there exists a radially symmetric global-in-time solution $(u,v)$ to \eqref{0} blowing up in infinite time with respect to the $L^\infty$-norm.
\end{theo}
\begin{rem}\label{Uwaga}
Notice that for $\alpha\in(\frac{2}{n},1)$ in \eqref{G1}, choosing $\phi(u)=\beta(u)$ we make sure that \eqref{G1} and
\eqref{H1} are satisfied mutually with \eqref{balance} and \eqref{decay} indicating that the assumptions of Theorem~\ref{theo1.5} are not contradictory.      
\end{rem}

\mysection{Preliminaries}\label{section2}
In this section we state some known results concerning local existence of solutions to \eqref{0} as well as some useful
properties of the solutions.
\begin{lem}\label{lem2.1}
  Suppose that $(u_0,v_0) \in C^0(\bar\Omega) \times W^{1,\infty}(\Omega)$ are radially symmetric and positive in $\bar\Omega$,
  and let $q\in (n,\infty)$.
  Then there exist $\tm \in (0,\infty]$ and a classical solution $(u,v)$ of \eqref{0} in
  $\Omega \times (0,\tm)$, where $u$ and $v$ are radially symmetric functions and satisfy
  \begin{eqnarray*}
	& & u \in C^0([0,\tm);C^0(\bar\Omega)) \cap C^{2,1}(\bar\Omega \times (0,\tm)), \\
	& & v \in C^0([0,\tm);W^{1,q}(\Omega)) \cap C^{2,1}(\bar\Omega \times (0,\tm)).
  \end{eqnarray*}
  Moreover, 
  \bas
	\mbox{either $\tm=\infty$,} \qquad \mbox{or } \|u(\cdot,t)\|_{L^\infty(\Omega)} 		
	\to \infty \quad \mbox{as } t\nearrow \tm
  \eas
  is fulfilled, equation \eqref{liapunov} holds and we have
  \begin{eqnarray}
	\io u(x,t)dx &=& \io u_0 \qquad \mbox{for all } t\in (0,\tm), \label{mass}\\
	\io v(x,t)dx &\le& \max \Big\{ \io u_0, \io v_0 \Big\}
	\qquad \mbox{for all } t\in (0,\tm). \label{mass_v}
  \end{eqnarray}
\end{lem}

\proof The claims concerning existence and regularity of the solution follow from well-known parabolic regularity theory and
  fixed point arguments, and the extensibility criterion also is proved by standard arguments. For details, we refer the 
  reader to \cite{amann, horstmann_winkler, wrzosek}. Moreover, the energy equation \eqref{liapunov} is proved in
  \cite[Lemma~2.1]{win_mmas} and the mass identities \eqref{mass} and \eqref{mass_v} immediately follow from integrating the
  first and second equation in \eqref{0}, respectively, by using the Neumann boundary conditions along with an ODE comparison. Conservation of radial symmetry is a consequence of uniqueness of solutions and the adequate form of 
equations in \eqref{0}.   
\qed

Next, we state a consequence of the Gagliardo-Nirenberg and the Young inequalities which will be used in forthcoming proofs
and which is given in \cite[Lemma~2.2]{win_bu} (see \cite{friedman} for details of the proof).
\begin{lem}\label{lem2.2}
  There is $C>0$ such that
  \be{2.2.1}
	\|\varphi\|_{L^2(\Omega)} \le C\|\nabla \varphi\|_{L^2(\Omega)}^\frac{n}{n+2} \|\varphi\|_{L^1(\Omega)}^\frac{2}{n+2}
	+ C\|\varphi\|_{L^1(\Omega)}
	\qquad \mbox{for all } \varphi \in W^{1,2}(\Omega).
  \ee
  In addition, for any $\eps>0$ there exists $C(\eps)>0$ such that
  \be{2.2.2}
	\|\varphi\|_{L^2(\Omega)}^2 \le \eps \|\nabla \varphi\|_{L^2(\Omega)}^2 + C(\eps) \|\varphi\|_{L^1(\Omega)}^2
	\qquad \mbox{for all } \varphi \in W^{1,2}(\Omega).
  \ee
\end{lem}
The following pointwise upper bound for the function $v$ will be an important ingredient to prove finite-time blowup.
 The result is given in \cite[Corollary~3.3]{win_bu} and its proof 
is exactly the same as the one performed in \cite[Section~3]{win_bu} since there only the second equation in \eqref{0} is used.
\begin{lem}\label{lem2.3}
  Let $p\in (1,\frac{n}{n-1})$. Then there is $C(p)>0$ such that
  whenever $u_0 \in C^0(\bar\Omega)$ and $v_0 \in W^{1,\infty}(\Omega)$ are positive in $\bar{\Omega}$ and radially symmetric, 
  the solution of \eqref{0} satisfies
  \be{2.3.1}
	v(r,t) \le C(p) \cdot \Big( \|u_0\|_{L^1(\Omega)} + \|v_0\|_{L^1(\Omega)}
	+ \|\nabla v_0\|_{L^2(\Omega)} \Big) \cdot r^{-\frac{n-p}{p}}
  \ee	
  for all $(r,t) \in (0,R) \times (0,\tm)$.
\end{lem}
\mysection{Finite-time blowup: estimates for the Liapunov functional}\label{section3}

In this section, we estimate the Liapunov functional $\F$ in terms of the dissipation rate $\D$ and frequently use the
ideas from \cite[Section~4]{win_bu}, where the case $\phi(u) =1$ and $\psi(u) =u$ is studied. In order to be able to handle the more general system \eqref{0}, we introduce new estimates in Lemma~\ref{lem3.4} along with a more
careful choice of some constants and the use of the terms contained in $\F$ which were not used in \cite{win_bu}. 

Following the ansatz of \cite{win_bu}, in view of the previous section we fix $m>0$, $M>0$, $B>0$, and $\kappa>n-2$ and 
assume that
\be{m}
	\io u = m \qquad \mbox{and} \qquad \io v \le M
\ee
and 
\be{B}
	v(x) \le B|x|^{-\kappa} \qquad \mbox{for all } x\in\Omega,
\ee
are satisfied. Moreover, we define the space
\bea{S}
	\set &:=& \bigg\{ (u,v) \in C^1(\bar\Omega) \times C^2(\bar\Omega) \ \bigg| \
	\mbox{$u$ and $v$ are positive and radially} \nn\\
	& & \hspace*{5mm}
	\mbox{symmetric satisfying $\frac{\partial v}{\partial\nu}=0$ on $\pO$, \eqref{m}, and \eqref{B}} \bigg\}.
\eea
The goal of this section is to prove that the inequality
\begin{equation}\label{4.1}
  \frac{\F(u,v)}{\D^\theta(u,v)+1} \ge - C(m,M,B,\kappa) \qquad \mbox{for all } (u,v) \in \set
\end{equation} 
holds with some constants $\theta \in (0,1)$ and $C(m,M,B, \kappa)>0$ (see Theorem~\ref{theo3.6}). Here it will be important
to state precisely the dependence of $C$ on $M$ and $B$. 

The main ingredient of the proof of \eqref{4.1} is the following estimate of $\io uv$.
\begin{lem}\label{lem3.1}
  Let \eqref{H1} and \eqref{psi} be fulfilled. Then there are $C(m,\kappa)>0$ and 
  \be{theta}
	\theta:=\frac{1}{1+\frac{n}{(2n+4)\kappa}} \, 
	\in \Big(\frac{1}{2},1\Big)
  \ee
  such that all $(u,v)\in\set$ satisfy
  \begin{eqnarray}\label{3.1.1}
	\io uv &\le&  C(m,\kappa) \cdot \left( 1+M^2 + B^{\frac{2n+4}{n+4}} \right) \cdot 
	\Bigg( \Big\|\Delta v-v+u\Big\|_{L^2(\Omega)}^{2\theta} \nn \\
	& & + \left\|\frac{\phi (u)}{\sqrt{\psi(u)}}\nabla u -\sqrt{\psi (u)}\nabla v\right\|_{L^2(\Omega)} +1 \Bigg) 
    +\frac{1}{2} \io | \nabla v |^2 + \io G(u).
  \end{eqnarray}
\end{lem} 
Lemma~\ref{lem3.1} is a generalization of \cite[Lemma~4.1]{win_bu} and our proof, which will be given after proving
several claims in the forthcoming lemmata, is based on the ideas given in 
\cite[Section~4]{win_bu} along with some additional estimates in order to cope with the more general functions $\phi$ and 
$\psi$. 

For notational convenience, we abbreviate
\be{f}
	f:=-\Delta v + v - u
\ee
and
\be{g}
	g:= \left( \frac{\phi (u)}{\sqrt{\psi(u)}}\nabla u -\sqrt{\psi (u)}\nabla v \right) \cdot \frac{x}{|x|}, \qquad x \neq 0,
\ee
for $(u,v)\in \set$.

The first step towards the proof of Lemma~\ref{lem3.1} is the following estimate which is completely similar to 
\cite[Lemma~4.2]{win_bu}. But as our different choice of the constants and their precise dependence on $M$ are important for 
the sequel, we give the proof for the reader's convenience.
\begin{lem}\label{lem3.2}
  For any $\eps \in (0,1)$ there exists $C(\eps)>0$ such that for all $(u,v) \in \set$ 
  \be{3.2.1}
	\io uv \le (1+ \eps) \io |\nabla v|^2 + C(\eps) \cdot \left( 1+ M^2 \right) \cdot \left( 
	\Big\| \Delta v - v + u \Big\|_{L^2(\Omega)}^{\frac{2n+4}{n+4}} +1 \right)
  \ee
  is fulfilled.
\end{lem}
\proof
  Multiplying \eqref{f} by $v$ and integrating by parts over $\Omega$ we have
  \be{3.2.2}
	\io uv = \io |\nabla v|^2 + \io v^2 - \io fv.
  \ee
  Now given $\eps \in (0,1)$, by Lemma \ref{lem2.2} and \eqref{m} we can fix 
  $c_1 = C_1 \cdot (1+M) >0$ and $c_2=C_2(\eps) \cdot M^2>0$ such that
  \be{3.2.3}
	\|v\|_{L^2(\Omega)} \le c_1 \cdot \left( \|\nabla v\|_{L^2(\Omega)}^\frac{n}{n+2}+1 \right)
  \ee
  and
  \be{3.2.4}
	\io v^2 \le \frac{\eps}{2} \io |\nabla v|^2 + c_2.
  \ee
  Applying the Cauchy-Schwarz inequality along with (\ref{3.2.3}) and Young's inequality
  (with exponents $\frac{2n+4}{n}$ and $\frac{2n+4}{n+4}$), we obtain $c_3=C_3(\eps) \cdot (1+M^{\frac{2n+4}{n+4}})>0$ 
  such that
  \bea{3.2.5}
	- \io fv
	&\le& \|f\|_{L^2(\Omega)} \|v\|_{L^2(\Omega)} 
	\le c_1 \cdot (\|\nabla v\|_{L^2(\Omega)}^\frac{n}{n+2}+1) \cdot \|f\|_{L^2(\Omega)} \nn\\
	&\le& \frac{\eps}{2} \io |\nabla v|^2 
	+ c_3 \|f\|_{L^2(\Omega)}^\frac{2n+4}{n+4} + c_1 \|f\|_{L^2(\Omega)}.
  \eea
  Since $\frac{2n+4}{n+4}>1$, we use Young's inequality once more and deduce that 
  \bas
	c_1 \|f\|_{L^2(\Omega)} \le c_1 \|f\|_{L^2(\Omega)}^\frac{2n+4}{n+4} + c_1
  \eas
  is satisfied. Combining the latter inequality with \eqref{3.2.2}, \eqref{3.2.4}, and \eqref{3.2.5},
  the claimed estimate \eqref{3.2.1} is proved, where we use $\frac{2n+4}{n+4} <2$ to deduce the estimate $(1+M^2)$ in
  \eqref{3.2.1}.
\qed

In view of Lemma~\ref{lem3.2}, the next step is to estimate $\int |\nabla v|^2$. This is first done in the annulus
$\Omega \setminus B_{r_0}$, where the value of $r_0$ will be fixed in Lemma~\ref{lem3.5} below. Since in \cite[Lemma~4.3]{win_bu} only equation \eqref{f} is used we could simply repeat its proof. However we give it in details
in order to state the exact dependence of the constants on $M$ and $B$ which will be of importance further.  
\begin{lem}\label{lem3.3}
  For any $r_0\in (0,R)$ and $\eps \in (0,1)$, there exists a constant 
  $C(\eps,m,\kappa)>0$ such that all $(u,v)\in\set$ satisfy
  \bea{3.3.1}
	\int_{\Omega \setminus B_{r_0}} |\nabla v|^2 
	&\le& \eps \io uv + \eps \io |\nabla v|^2 
    + C(\eps,m,\kappa) \cdot \left( 1+ M^{\frac{2n+4}{n+4}} + B^{\frac{2n+4}{n+4}} \right) \cdot \Bigg\{
	r_0^{-\frac{2n+4}{n}\kappa} \nn \\
	& & + \Big\|\Delta v-v+u\Big\|_{L^2(\Omega)}^\frac{2n+4}{n+4} \Bigg\}.
  \eea
\end{lem}
\proof
  Let $\alpha_1 \in (0,1)$ be arbitrary.
  As $v > 0$, a multiplication of \eqref{f} by $v^{\alpha_1}$ and an integration by parts
  over $\Omega$ implies 
  \be{3.3.2}
	\alpha_1 \io v^{\alpha_1-1}|\nabla v|^2 \le \alpha_1 \io v^{\alpha_1-1}|\nabla v|^2 + \io v^{\alpha_1+1}
	= \io uv^{\alpha_1} + \io fv^{\alpha_1}.
  \ee
  Using next \eqref{B} and $\alpha_1 \in (0,1)$, we obtain
  \bas
	\alpha_1 \io v^{\alpha_1-1} |\nabla v|^2 \ge \alpha_1 B^{\alpha_1-1} r_0^{(1-\alpha_1)\kappa} \cdot
	\int_{\Omega \setminus B_{r_0}} |\nabla v|^2, 
  \eas
  whence \eqref{3.3.2} yields
  \bea{3.3.3}
	\int_{\Omega \setminus B_{r_0}} |\nabla v|^2 
	&\le& \frac{B^{1-\alpha_1}}{\alpha_1} r_0^{-(1-\alpha_1)\kappa} \io uv^{\alpha_1} \,
	+ \, \frac{B^{1-\alpha_1}}{\alpha_1} r_0^{-(1-\alpha_1)\kappa} \io fv^{\alpha_1}.
  \eea
  In view of $\alpha_1 \in (0,1)$ and Young's inequality, for any $\eta>0$ there is $c_1(\eta,B) = C_1 (\eta) \cdot B>0$ such 
  that
  \be{3.3.4}
	\frac{B^{1-\alpha_1}}{\alpha_1} r_0^{-(1-\alpha_1)\kappa} v^{\alpha_1}(r)
	\le \eta v(r) + c_1(\eta,B) r_0^{-\kappa}
	\qquad \mbox{for all } r\in (0,R).
  \ee
  The choice $\eta:=\eps$ implies
  \bea{3.3.5}
	\frac{B^{1-\alpha_1}}{\alpha_1} r_0^{-(1-\alpha_1)\kappa} \io uv^{\alpha_1}
	&\le& \eps \io uv 
	+ c_1(\eps,B) r_0^{-\kappa} \io u \nn\\
	 &=& \eps \io uv 
	+ c_1(\eps,B) m r_0^{-\kappa} \nn\\
	&\le& \eps \io uv 
	+ c_1(\eps,B) m R^{\frac{n+4}{n}\kappa} r_0^{-\frac{2n+4}{n}\kappa}
  \eea
  due to \eqref{m} and $u \ge 0$.
  
  Furthermore, using (\ref{3.3.4}) with $\eta:=1$ along with the Cauchy-Schwarz inequality, we deduce that
  \begin{eqnarray}\label{3.3.6}
	\frac{B^{1-\alpha_1}}{\alpha_1} r_0^{-(1-\alpha_1)\kappa} \io fv^{\alpha_1}
	&\le& \io |f|v + c_1(1,B) r_0^{-\kappa} \io |f|, \nn \\
	&\le& \|f\|_{L^2(\Omega)} \|v\|_{L^2(\Omega)} + c_1(1,B) r_0^{-\kappa} \sqrt{|\Omega|} \|f\|_{L^2(\Omega)}.
  \end{eqnarray}
  Since by Lemma \ref{lem2.2} and \eqref{m}, there exists $c_2(M) = C_2 \cdot (1+M) >0$ such that
  \bas
	\|v\|_{L^2(\Omega)} \le c_2(M) \cdot \Big( \|\nabla v\|_{L^2(\Omega)}^\frac{n}{n+2} + 1 \Big)
    \le c_2(M) \cdot \Big( \|\nabla v\|_{L^2(\Omega)}^\frac{n}{n+2} + R^\kappa r_0^{-\kappa} \Big),
  \eas
  from \eqref{3.3.6} we infer
  \bas
	\frac{B^{1-\alpha_1}}{\alpha_1} r_0^{-(1-\alpha_1)\kappa} \io fv^{\alpha_1}
	\le c_3(M,B,\kappa) \cdot \Big(\|f\|_{L^2(\Omega)} \|\nabla v\|_{L^2(\Omega)}^\frac{n}{n+2}
	+r_0^{-\kappa} \|f\|_{L^2(\Omega)} \Big)
  \eas
  with some $c_3(M,B,\kappa) = C_3 (\kappa) \cdot (1+M+B)>0$.
  Applying Young's inequality, 
  \bas
	c_3(M,B,\kappa) \|f\|_{L^2(\Omega)} \|\nabla v\|_{L^2(\Omega)}^\frac{n}{n+2}
	\le \eps \|\nabla v\|_{L^2(\Omega)}^2
	+ c_4(\eps,M,B,\kappa) \|f\|_{L^2(\Omega)}^\frac{2n+4}{n+4}
  \eas
  and
  \bas
	c_3(M,B,\kappa) r_0^{-\kappa} \|f\|_{L^2(\Omega)}
	\le c_3(M,B,\kappa) \|f\|_{L^2(\Omega)}^\frac{2n+4}{n+4}
	+ c_3(M,B,\kappa) r_0^{-\frac{2n+4}{n}\kappa}
  \eas
  hold with some $c_4(\eps,M,B,\kappa) = C_4 (\eps, \kappa) \cdot \big( 1+ M^{\frac{2n+4}{n+4}} + B^{\frac{2n+4}{n+4}} \big) 
  >0$.
  Thus, \eqref{3.3.6} finally turns into
  \bas
	\frac{B^{1-\alpha_1}}{\alpha_1} r_0^{-(1-\alpha_1)\kappa} \io fv^{\alpha_1}
	&\le& \eps \io |\nabla v|^2
	+ \Big( c_4(\eps,M,B,\kappa) +c_3(M,B,\kappa) \Big) \cdot \|f\|_{L^2(\Omega)}^\frac{2n+4}{n+4} \\
	& & + c_3(M,B,\kappa) r_0^{-\frac{2n+4}{n}\kappa}.
  \eas
  In conjunction with \eqref{3.3.3} and \eqref{3.3.5}, the claim \eqref{3.3.1} is proved.
\qed

Next we prove a corresponding estimate of $\nabla v$ on the ball $B_{r_0}$. Our proof is based on ideas from 
\cite[Lemma~4.4]{win_bu} which are generalized to the problem \eqref{0}. We recall that $G$ and $H$ are defined in
\eqref{GH} and remark that the following proof is the only place where we use the assumption \eqref{psi}. Moreover, it is
important that $r_0$ can be chosen arbitrarily small in order to obtain a subquadratic power of $\| f \|_{L^2(\Omega)}$
in Lemma~\ref{lem3.5}. 
\begin{lem}\label{lem3.4}
  Assume that \eqref{H1} and \eqref{psi} are satisfied. Then there exist $\mu = \mu(\gamma) \in (0,2)$ 
  and $C(m)>0$ such that for all $r_0\in (0,R)$ and $(u,v)\in\set$ 
  \begin{eqnarray}\label{3.4.1}
	\int_{B_{r_0}} |\nabla v|^2
	&\le& \mu \io G(u) + C(m) \cdot \Bigg\{
	r_0 \cdot \Big\|\Delta v-v+u \Big\|_{L^2(\Omega)}^2 \nn \\
	& & +\left\| \frac{\phi (u)}{\sqrt{\psi(u)}}\nabla u -\sqrt{\psi (u)}\nabla v \right\|_{L^2(\Omega)} 
	+\|v\|_{L^2(\Omega)}^2
	+1 \Bigg\}
  \end{eqnarray}
  is fulfilled.
\end{lem} 
\proof As \eqref{H1} implies $(\frac{4(n-1)}{n-2} -2)\gamma <2$, we fix $\delta \in (0, \frac{2n-2}{R}]$ small 
enough such that 
\begin{equation}\label{mu}
 \mu := \left( \frac{4(n-1)}{n-2} e^{\delta R} -2 \right) \cdot \gamma  \in (0,2)
\end{equation}
is fulfilled. 
As $u$ and $v$ are radially symmetric, \eqref{f} and \eqref{g} 
yield the identities 
\begin{equation}\label{pier}
(r^{n-1} v_r)_r = -r^{n-1}u - r^{n-1}f + r^{n-1}v
\end{equation}
and 
\begin{equation}\label{pierd}
v_r = \frac{\phi(u)}{\psi(u)} u_r- \frac{g}{\sqrt{\psi(u)}}\;.
\end{equation} 
Multiplying (\ref{pier}) by $r^{n-1}v_r$ and using (\ref{pierd}) as well 
as Young's inequality, we obtain
\begin{eqnarray}\label{3.4.2}
 \frac{1}{2} \left( (r^{n-1} v_r)^2 \right)_r &=& -r^{2n-2}uv_r - r^{2n-2} fv_r + r^{2n-2} vv_r \nn \\
 &\le& -r^{2n-2} \frac{u \phi(u)}{\psi(u)} u_r + r^{2n-2} \frac{u}{\sqrt{\psi(u)}} g 
 + \frac{\delta}{2} (r^{n-1}v_r)^2 + \frac{1}{2 \delta} r^{2n-2} f^2 \nn \\
 & & + \frac{1}{2} r^{2n-2} (v^2)_r \qquad \mbox{for all } r \in (0,R).
\end{eqnarray}
Defining $y(r) := (r^{n-1}v_r)^2$, $r \in [0,R]$, we obtain 
$$y_r \le -2r^{2n-2} \frac{u \phi(u)}{\psi(u)} u_r + 2r^{2n-2} \frac{u}{\sqrt{\psi(u)}} g 
 + \delta y + \frac{1}{\delta} r^{2n-2} f^2 + r^{2n-2} (v^2)_r, \qquad r \in (0,R),$$ 
along with $y(0) =0$ due to the regularity of $v$. Thus, an integration implies
\begin{eqnarray}\label{3.4.3}
 r^{2n-2} v_r^2(r) = y(r) &\le&
	-2\int_0^r e^{\delta (r-\rho)} \rho^{2n-2} \frac{u(\rho) \phi(u(\rho))}{\psi(u(\rho))} u_r (\rho) 
        \; d\rho \nn \\ & &
	+2\int_0^r e^{\delta (r-\rho)} \rho^{2n-2} \frac{u(\rho)}{\sqrt{\psi(u(\rho))}} g (\rho) \; d\rho \nn\\
	& & +\frac{1}{\delta} \int_0^r e^{\delta (r-\rho)} \rho^{2n-2} f^2(\rho) d\rho
	+\int_0^r e^{\delta (r-\rho)} \rho^{2n-2} (v^2)_r(\rho) d\rho
\end{eqnarray}
for all $r \in (0,R)$.
Integrating by parts and using the nonnegativity of $H$ (defined in \eqref{GH}), we obtain
\begin{eqnarray}\label{3.4.4}
  && \hspace*{-20mm} 
  -2\int_0^r e^{\delta (r-\rho)} \rho^{2n-2} \frac{u (\rho) \phi(u (\rho))}{\psi(u (\rho))} u_r (\rho) \; d\rho 
  \nn \\
  &=& 4(n-1) \int_0^r e^{\delta (r-\rho)} \rho^{2n-3} H(u(\rho))  \; d\rho \nn \\
  & & -2 \delta \int_0^r e^{\delta (r-\rho)} \rho^{2n-2} H(u(\rho)) \; d\rho
  -2r^{2n-2} H(u(r)) \nn \\
  &\le& 4(n-1) e^{\delta R} \int_0^r \rho^{2n-3} H(u(\rho))  \; d\rho -2r^{2n-2} H(u(r)), \qquad r \in (0,R).
\end{eqnarray}
Next, denoting by $\omega_n$ the $(n-1)$-dimensional measure of the sphere $\partial B_1$ and applying the 
Cauchy-Schwarz inequality as well as \eqref{psi}, we deduce that 
\begin{eqnarray}\label{3.4.5}
 && \hspace*{-20mm} 
 2\int_0^r e^{\delta (r-\rho)} \rho^{2n-2} \frac{u(\rho)}{\sqrt{\psi(u(\rho))}} g (\rho) \; d\rho \nn \\
 &\le& 2 \left( \int_0^R \rho^{n-1} \frac{u^2(\rho)}{\psi(u(\rho))} \; d\rho \right)^\frac{1}{2} \cdot 
	\left( \int_0^r e^{2\delta(r-\rho)} \cdot \rho^{3n-3} g^2(\rho) \; d\rho \right)^\frac{1}{2} \nn \\
 &\le& 2 \left( \frac{1}{c_0} \int_0^R \rho^{n-1} u(\rho) \; d\rho \right)^\frac{1}{2} \cdot 
	\left( e^{2\delta R} r^{2n-2} \int_0^R \rho^{n-1} g^2(\rho) \; d\rho \right)^\frac{1}{2} \nn \\ 
 &\le& \frac{2 e^{\delta R}}{w_n \sqrt{c_0}} \cdot \sqrt{m} \cdot r^{n-1} \cdot \| g \|_{L^2(\Omega)},
 \qquad r \in (0,R).
\end{eqnarray}
Similarly, we estimate the third term on the right-hand side of \eqref{3.4.3} according to 
\begin{eqnarray}\label{3.4.6}
  \frac{1}{\delta} \int_0^r e^{\delta (r-\rho)} \rho^{2n-2} f^2(\rho) \; d\rho
	&\le& \frac{e^{\delta R}}{\delta} \cdot r^{n-1} \cdot \int_0^R \rho^{n-1} f^2(\rho) \; d\rho \nn\\
	&=& \frac{e^{\delta R}}{\delta \omega_n} \cdot r^{n-1} \cdot \|f\|_{L^2(\Omega)}^2
	\qquad \mbox{for all } r\in (0,R).
\end{eqnarray}
As $\delta \le \frac{2n-2}{R}$ yields $(2n-2) \rho^{2n-3} \ge \delta \rho^{2n-2}$ for all $\rho \in (0,R)$,
integrating by parts we furthermore arrive at 
\begin{eqnarray}\label{3.4.7}
 \int_0^r e^{\delta(r-\rho)} \rho^{2n-2} (v^2)_r(\rho) \; d\rho
	&=& r^{2n-2} v^2(r) \nn\\
	& & - \int_0^r e^{\delta(r-\rho)} \cdot [ (2n-2) \rho^{2n-3} - \delta \rho^{2n-2}] \cdot v^2(\rho) \;
        d\rho \nn\\
	&\le& r^{2n-2} v^2(r)
	\qquad \mbox{for all } r\in (0,R).
\end{eqnarray}
Hence, \eqref{3.4.3}-\eqref{3.4.7} imply that there is a constant $c_1(m) >0$ such that
\begin{eqnarray*}
  r^{2n-2} v_r^2(r) &\le& 4(n-1) e^{\delta R} \int_0^r \rho^{2n-3} H(u(\rho))  \; d\rho -2r^{2n-2} H(u(r)) \\
  & & + \frac{c_1(m)}{\omega_n} r^{n-1} \|g\|_{L^2(\Omega)} 
   + \frac{c_1(m)}{\omega_n} r^{n-1} \|f\|_{L^2(\Omega)}^2
	+ r^{2n-2} v^2(r), \qquad r \in (0,R).
\end{eqnarray*}
Multiplying this inequality by $\omega_n r^{1-n}$ and integrating over $r \in (0,r_0)$, we have
\begin{eqnarray}\label{3.4.8}
 \int_{B_{r_0}} \| \nabla v \|^2 &=& \omega_n \int_0^{r_0} r^{n-1} v_r^2(r) \; dr \nn \\
 &\le& 4(n-1) e^{\delta R} \omega_n  \int_0^{r_0} r^{1-n} \int_0^r \rho^{2n-3} H(u(\rho))  \; d\rho \, dr \nn \\
 & &  -2\omega_n \int_0^{r_0} r^{n-1} H(u(r)) \; dr + c_1 (m) r_0 \|g\|_{L^2(\Omega)}  \nn \\
 & &  + c_1(m) r_0 \|f\|_{L^2(\Omega)}^2 + \omega_n \int_0^{r_0} r^{n-1} v^2(r) \; dr \nn \\
 &\le& 4(n-1) e^{\delta R} \omega_n  \int_0^{r_0} r^{1-n} \int_0^r \rho^{2n-3} H(u(\rho))  \; d\rho \, dr \nn \\
 & &  -2 \int_{B_{r_0}} H(u) + c_1 (m) R \|g\|_{L^2(\Omega)}  
 + c_1(m) r_0 \|f\|_{L^2(\Omega)}^2 + \|v\|_{L^2(\Omega)}^2.
\end{eqnarray}
Finally, Fubini's theorem, $n \ge 3$, the nonnegativity of $H$, \eqref{H1}, and \eqref{mu} yield 
\begin{eqnarray*}
 && \hspace*{-20mm}
   4(n-1) e^{\delta R} \omega_n  \int_0^{r_0} r^{1-n} \int_0^r \rho^{2n-3} H(u(\rho))  \; d\rho \, dr 
   -2 \int_{B_{r_0}} H(u) \\
 &=& 4(n-1) e^{\delta R} \omega_n  \int_0^{r_0} \left( \int_{\rho}^{r_0} r^{1-n} \; dr \right) 
    \rho^{2n-3} H(u(\rho))  \; d\rho -2 \int_{B_{r_0}} H(u) \\
 &=& \frac{4(n-1)}{n-2} e^{\delta R} \omega_n  \int_0^{r_0} \left( \rho^{2-n} - r_0^{2-n} \right) 
    \rho^{2n-3} H(u(\rho))  \; d\rho -2 \int_{B_{r_0}} H(u) \\
 &\le& \frac{4(n-1)}{n-2} e^{\delta R} \omega_n  \int_0^{r_0}  
    \rho^{n-1} H(u(\rho))  \; d\rho -2 \int_{B_{r_0}} H(u) \\
 &=& \left( \frac{4(n-1)}{n-2} e^{\delta R} -2 \right) \int_{B_{r_0}} H(u) 
 \le \left( \frac{4(n-1)}{n-2} e^{\delta R} -2 \right) \int_{\Omega} H(u) \\
 &\le& \left( \frac{4(n-1)}{n-2} e^{\delta R} -2 \right) \int_{\Omega} \left( \gamma G(u) + b(u+1) \right) 
 = \mu \io G(u) + c_2 (m) 
\end{eqnarray*}
with some $c_2(m) >0$. Upon a combination with \eqref{3.4.8}, the claim is proved.
\qed

The final step towards the proof of \eqref{3.1.1} is now a combination of Lemma~\ref{lem3.3} and 
Lemma~\ref{lem3.4}. The proof is very similar to the one given in \cite[Lemma~4.5]{win_bu}, but as we have to 
choose some constants in a different way, we give the proof for completeness of our arguments.
\begin{lem}\label{lem3.5}
  Suppose that \eqref{H1} and \eqref{psi} are fulfilled and let $\theta \in (\frac{1}{2},1)$ and 
  $\mu \in (0,2)$ be as defined in \eqref{theta} and \eqref{mu},
  respectively. Then for any $\eps \in (0, \frac{1}{2})$ there exists $C(\eps,m,\kappa)>0$ such that 
  \begin{eqnarray}\label{3.5.1}
	\io |\nabla v|^2 
	&\le&  C(\eps,m,\kappa) \cdot \left( 1+ M^2 
    + B^{\frac{2n+4}{n+4}} \right) \cdot 
	\bigg( \Big\|\Delta v-v+u\Big\|_{L^2(\Omega)}^{2\theta} \nn \\ & &
	+ \Big\|\frac{\phi (u)}{\sqrt{\psi(u)}}\nabla u -\sqrt{\psi (u)}\nabla v\Big\|_{L^2(\Omega)} +1 \bigg) \nn \\ &&
	+ \frac{\eps}{1-2\eps} \io uv + \frac{\mu}{1-2\eps} \io G(u)
  \end{eqnarray}
  is fulfilled for all $(u,v) \in \set$.
\end{lem} 
\proof We fix $\eps \in (0, \frac{1}{2})$ and set $\beta := \frac{(2n+4)\kappa}{n}$ which implies 
$\theta = \frac{\beta}{\beta+1}$. Next we define $r_0 := 
\min\{\frac{R}{2}, \|f\|_{L^2(\Omega)}^{-\frac{2}{\beta+1}} \} \in (0,R)$. Hence, by Lemma~\ref{lem3.3} there 
is $c_1= C_1 (\eps, m, \kappa) \cdot \big( 1+ M^{\frac{2n+4}{n+4}} + B^{\frac{2n+4}{n+4}} \big)>0$ such that
\begin{equation}\label{3.5.2}
 \int_{\Omega \setminus B_{r_0}} |\nabla v|^2
	\le \eps \io uv + \eps \io |\nabla v|^2
	+ c_1 \cdot \Big( r_0^{-\beta} + \|f\|_{L^2(\Omega)}^\frac{2n+4}{n+4} \Big).
\end{equation}
Applying next Lemma~\ref{lem3.4}, we get a constant $c_2 = c_2(m)$ such that
\begin{equation}\label{3.5.3}
 \int_{B_{r_0}} |\nabla v|^2  \le \mu \io G(u) + c_2 \cdot \Big( r_0 \|f\|_{L^2(\Omega)}^2 + \|g\|_{L^2(\Omega)} 
   + \|v\|_{L^2(\Omega)}^2 + 1 \Big).
\end{equation}
Adding both inequalities, we deduce that
\begin{eqnarray}\label{3.5.4}
 (1-\eps) \io |\nabla v|^2  &\le& \eps \io uv + \mu \io G(u) + c_1 r_0^{-\beta} + 
  c_1 \|f\|_{L^2(\Omega)}^\frac{2n+4}{n+4}  + c_2 r_0 \|f\|_{L^2(\Omega)}^2 \nn \\
  & & + c_2 (\|g\|_{L^2(\Omega)} +1) + c_2 \|v\|_{L^2(\Omega)}^2. 
\end{eqnarray}
Next, by Lemma~\ref{lem2.2} and \eqref{m} there exists $c_3 = C_3 (\eps,m) \cdot M^2 >0$ such that
$$c_2 \|v\|_{L^2(\Omega)}^2 \le \eps \io |\nabla v|^2 + c_3,$$ 
which inserted into \eqref{3.5.4} yields 
\begin{eqnarray}\label{3.5.5}
 (1-2\eps) \io |\nabla v|^2  &\le& \eps \io uv + \mu \io G(u) + c_2 (\|g\|_{L^2(\Omega)} +1) + c_3  +I,
\end{eqnarray}
where we set 
$$I := c_1 r_0^{-\beta} + 
  c_1 \|f\|_{L^2(\Omega)}^\frac{2n+4}{n+4}  + c_2 r_0 \|f\|_{L^2(\Omega)}^2.$$
In case of $\|f\|_{L^2(\Omega)} \le (\frac{2}{R})^\frac{\beta+1}{2}$, we have $r_0=\frac{R}{2}$ and conclude that
  \bas
	I \le c_1 \cdot \Big(\frac{2}{R}\Big)^\beta
	+ c_1 \cdot \Big(\frac{2}{R}\Big)^{\frac{\beta+1}{2} \cdot \frac{2n+4}{n+4}}
	+ c_2 \cdot \frac{R}{2} \cdot \Big( \frac{2}{R} \Big)^{\beta+1},
  \eas
  which in conjunction with \eqref{3.5.5} proves \eqref{3.5.1} in this case.

Furthermore, in the case $\|f\|_{L^2(\Omega)} > (\frac{2}{R})^\frac{\beta+1}{2}$ we have 
$r_0=\|f\|_{L^2(\Omega)}^{-\frac{2}{\beta+1}}$ and therefore  
$$ I \le c_1 \|f\|_{L^2(\Omega)}^\frac{2\beta}{\beta+1}
	+c_1 \|f\|_{L^2(\Omega)}^\frac{2n+4}{n+4}
	+ c_2 \|f\|_{L^2(\Omega)}^{2-\frac{2}{\beta+1}} 
	= (c_1+c_2) \|f\|_{L^2(\Omega)}^\frac{2\beta}{\beta+1} 
	+c_1 \|f\|_{L^2(\Omega)}^\frac{2n+4}{n+4}.$$
In view of $\kappa > n-2$ and $n \ge 3$, we calculate
$$\frac{\beta}{\frac{n+2}{2}} =\frac{2}{n+2} \cdot \frac{(2n+4)\kappa}{n}
	>\frac{4(n-2)}{n} \ge \frac{4}{3}>1$$
which implies that $2\theta = \frac{2\beta}{\beta+1}>\frac{2n+4}{n+4}$. Applying once more Young's inequality, 
we obtain 
$$I \le (2c_1+c_2) \|f\|_{L^2(\Omega)}^\frac{2\beta}{\beta+1} +  c_1,$$
which inserted into \eqref{3.5.5} proves \eqref{3.5.1} in the case 
  $\|f\|_{L^2(\Omega)} > (\frac{2}{R})^\frac{\beta+1}{2}$ and thereby completes the proof.
\qed

Next, we complete the proof of the announced estimate \eqref{3.1.1}. \\

{\bf Proof of Lemma~\ref{lem3.1}.} \quad Let $\mu \in (0,2)$ be as defined in \eqref{mu}. In view of $\mu <2$
there exists $\eta \in (0, \frac{1}{2})$ such that $\mu(1-\eta) <1$. Keeping this value of $\eta$ fixed, we 
moreover fix $\eps \in (0, \frac{1}{4})$ small enough such that
\begin{equation}\label{3.1.2}
 \frac{\mu (1+ \eps - \eta)}{1- 3 \eps -\eps^2 + \eps \eta} \le 1
 \quad\mbox{and}\qquad \frac{\eta(1-2\eps)}{1- 3 \eps -\eps^2 + \eps \eta} \le \frac{1}{2}.
\end{equation}
An application of Lemma~\ref{lem3.2} implies the existence of $c_1 = C_1 \cdot (1+M^2)>0$ such that
$$\io uv \le \eta \io |\nabla v|^2  + (1+ \eps - \eta) \io |\nabla v|^2 + c_1 \cdot \left( 
	\| f \|_{L^2(\Omega)}^{\frac{2n+4}{n+4}} +1 \right).$$
Furthermore, by Lemma~\ref{lem3.5} there is $c_2 = C_2(m,\kappa) \cdot \big( 1+ M^2 
    + B^{\frac{2n+4}{n+4}} \big)>0$ such that
\begin{eqnarray*}
\io uv &\le& \eta \io |\nabla v|^2  +  \frac{\eps(1+ \eps - \eta)}{1-2\eps} \io uv 
  + \frac{\mu(1+ \eps - \eta)}{1-2\eps} \io G(u) \nn \\ & & + c_2(1+ \eps - \eta) \cdot 
	\left(\|f\|_{L^2(\Omega)}^{2\theta} 
	+ \|g\|_{L^2(\Omega)} +1 \right) + c_1 \cdot \left( 
	\| f \|_{L^2(\Omega)}^{\frac{2n+4}{n+4}} +1 \right).
\end{eqnarray*}
A rearrangement of the terms yields 
\begin{eqnarray*}
\io uv &\le& \frac{\eta(1-2\eps)}{1- 3 \eps -\eps^2 + \eps \eta} \io |\nabla v|^2   
  + \frac{\mu(1+ \eps - \eta)}{1- 3 \eps -\eps^2 + \eps \eta} \io G(u) \nn \\ & & 
  + c_3 \cdot 
	\left(\|f\|_{L^2(\Omega)}^{2\theta} + \| f \|_{L^2(\Omega)}^{\frac{2n+4}{n+4}}
	+ \|g\|_{L^2(\Omega)} +1 \right) 
\end{eqnarray*}
with some $c_3 = C_3(m,\kappa) \cdot \big( 1+ M^2 
    + B^{\frac{2n+4}{n+4}} \big)>0$. As $\frac{2n+4}{n+4}<2\theta$ (which is shown in Lemma~\ref{lem3.5}),
a further application of the Young inequality along with \eqref{3.1.2} implies \eqref{3.1.1}.
\qed

The final result of this section is to show that the Liapunov functional $\F$ can be estimated according to
\eqref{4.1}. The proof uses the idea of \cite[Theorem~5.1]{win_bu} as a basic ingredient, but in fact our 
estimates also make use of the other terms which are contained in $\F$.
\begin{theo}\label{theo3.6}
  Assume that \eqref{H1} and \eqref{psi} are satisfied and let $\theta \in (\frac{1}{2},1)$ be as defined 
  in \eqref{theta}. Then there exists $C(m,\kappa)>0$ such that 
  \be{3.6.1}
	\F(u,v) \ge -C(m,\kappa) \cdot \left( 1+M^2 + B^{\frac{2n+4}{n+4}} \right) \cdot \Big( \D^\theta (u,v) +1 \Big)
  \ee
  is fulfilled for all $(u,v) \in \set$, where $\F$ and $\D$ are given in \eqref{F} and \eqref{D}, 
  respectively.
\end{theo}
\proof
  In view of \eqref{f}, \eqref{g}, and $\theta > \frac{1}{2}$, an application of Young's inequality
  to \eqref{3.1.1} implies the existence of $c_1=C_1(m,\kappa) \cdot \big( 1+M^2 + B^{\frac{2n+4}{n+4}} \big)>0$ such that
  \bas
	\io uv \le c_1 \left( \Big(\|f\|_{L^2(\Omega)}^2 + \|g\|_{L^2(\Omega)}^2 \Big)^\theta +1 \right)
        +\frac{1}{2} \io | \nabla v |^2 + \io G(u).
  \eas
  Hence, we conclude that
  \bas
	\F(u,v)
	&=& \frac{1}{2} \io |\nabla v|^2 + \frac{1}{2} \io v^2 - \io uv + \io G(u) \\
	&\ge& - c_1 \cdot \left( (\|f\|_{L^2(\Omega)}^2 + \|g\|_{L^2(\Omega)}^2)^\theta +1 \right).
  \eas
  As \eqref{D}, \eqref{f}, and \eqref{g} imply $\D(u,v)=\|f\|_{L^2(\Omega)}^2 + \|g\|_{L^2(\Omega)}^2$,
  the proof is complete.
\qed
\mysection{Finite-time blowup: proof of the main results}\label{section4}
In view of Theorem~\ref{theo3.6} and $\theta \in (0,1)$, we derive an ODI for the function $y(t) :=
-\F(u(\cdot,t),v(\cdot,t))$ with superlinear nonlinearity. This shows that the solution $(u,v)$ 
blows up in finite time if $-\F(u_0,v_0)$ is large. The following result and its proof 
are completely the same as in \cite[Lemma~5.2]{win_bu} so that we confine ourselves to giving only a sketch 
of the main ideas of the proof.
\begin{lem}\label{lem4.1}
  Suppose that \eqref{H1} and \eqref{psi} are fulfilled, let $\theta \in (\frac{1}{2},1)$ be as defined 
  in \eqref{theta}and let $m>0, A>0$ and $\kappa>n-2$. Then there exist $K=K(m,A,\kappa) = k(m,\kappa) \cdot (1+A^2)>0$ and 
  $C=C(m,A,\kappa)>0$ such that for any 
  \bea{tB}
    (u_0,v_0) \in \widetilde{\B}(m,A,\kappa) &:=& \bigg\{
	(u_0,v_0) \in C^0(\bar\Omega) \times W^{1,\infty}(\Omega) \ \bigg| \
	\mbox{$u_0$ and $v_0$ are radially symmetric} \nn \\ 
    & & \hspace*{5mm} \mbox{and positive in $\bar\Omega$, 
	$\io u_0=m$, $\|v_0\|_{W^{1,2}(\Omega)} \le A$,}\nn \\ 
    & & \hspace*{5mm} \mbox{and $\F(u_0,v_0) \le -K $} \bigg\}
  \eea
  the corresponding solution $(u,v)$ of \eqref{0} satisfies
  \be{4.1.1}
	\F(u(\cdot,t),v(\cdot,t))
	\le \frac{\F (u_0,v_0)}{(1-Ct)^\frac{\theta}{1-\theta}}
	\qquad \mbox{for all } t\in (0,\tm).
  \ee
  In particular, $(u,v)$ blows up in finite time $\tm \le \frac{1}{C}$.
\end{lem}
\proof We only give a sketch of the main ideas and refer to \cite[Lemma~5.2]{win_bu} for further details.

We fix $c_1 >0$ such that
$$\|\varphi\|_{L^1(\Omega)} \le c_1 \|\varphi\|_{W^{1,2}(\Omega)}
	\qquad \mbox{for all } \varphi \in W^{1,2}(\Omega).$$
Moreover, in view of $\kappa > n-2$ and Lemma~\ref{lem2.3}, there is $c_2 = c_2 (\kappa)>0$ such that 
for any $(u_0,v_0) \in \widetilde{\B}(m,A,\kappa)$ the solution $(u,v)$ to \eqref{0} fulfills
\begin{equation}\label{4.1.2}
  v(r,t) \le c_2 \cdot \Big( \|u_0\|_{L^1(\Omega)} + \|v_0\|_{L^1(\Omega)} + \|\nabla v_0\|_{L^2(\Omega)} \Big)
	\cdot r^{-\kappa}
\end{equation}
for all $(r,t) \in (0,R) \times (0,\tm)$. Setting $B:=c_2 (m+c_1 A+A)$ and $M:=\max \{m,c_1 A\}$, 
Lemma~\ref{lem2.1} and \eqref{4.1.2} imply that $(u(\cdot,t),v(\cdot,t)) \in \set$
  for all $t\in (0,\tm)$ provided that $(u_0,v_0) \in \widetilde{\B}(m,A,\kappa)$. 
In view of Theorem~\ref{theo3.6} and our definition of $B$ and $M$, there is a constant $c_3 = C_3(m,\kappa) \cdot \left( 1+
A^2 \right)$ such that 
\begin{equation}\label{4.1.3}
  \F(u(\cdot,t),v(\cdot,t)) \ge -c_3 \cdot \Big( \D^\theta (u(\cdot,t),v(\cdot,t)) +1 \Big)
\end{equation}  
is satisfied for all $t\in (0,\tm)$ provided that $(u_0,v_0) \in \widetilde{\B}(m,A,\kappa)$.
Hence, we set  
$K(m,A,\kappa) = 2c_3$, $C(m,A,\kappa) = \frac{1-\theta}{2c_3 \theta}$, and 
$y(t):=-\F(u(\cdot,t),v(\cdot,t))$, $t\in [0,\tm)$, for $(u_0,v_0) \in \widetilde{\B}(m,A,\kappa)$. 
As $y$ is nonincreasing by \eqref{liapunov} and therefore satisfies $y(t) \ge 2c_3$ for $t \in (0,\tm)$, 
\eqref{4.1.3} and \eqref{liapunov} imply 
$$y^\prime (t) \ge \Big( \frac{y(t)}{2c_3} \Big)^\frac{1}{\theta}
	\qquad \mbox{for all } t\in (0,\tm),$$
which implies \eqref{4.1.1}.
\qed

The proof of Theorem~\ref{theo1} is now immediate.\\

{\bf Proof of Theorem~\ref{theo1}.} \quad We fix an arbitrary $\kappa > n-2$. Then the claim directly follows
from Lemma~\ref{lem4.1} by defining $K(m,A) :=k(m,\kappa)$ and $T(m,A):=\frac{1}{C(m,A,\kappa)}$, where 
$k(m,\kappa)$ and $C(m,A,\kappa)$ are provided in Lemma~\ref{lem4.1}.
\qed

Let us next show that the set $\B(m,A)$ defined in \eqref{t1.1} has the properties claimed in Theorem~\ref{theo1.2}. 
Since the condition 
\begin{equation}\label{1.2.1} 
  \F(u_0,v_0) \le -K(m) \cdot (1+A^\tau)
\end{equation} 
in \eqref{t1.1} is given with $\tau =2$, we can use the functions constructed in \cite[Lemma~4.1]{win_mmas} to deduce that $\B(m,A) \neq \emptyset$ without any additional
restriction on $\alpha$ (which is given in \eqref{G1}). In case of $\tau >2$, this is not possible. Moreover, as 
\eqref{1.2.1} cannot be imposed for $\tau<2$ in view of the Liapunov functional $\F$, the condition \eqref{1.2.1} with
$\tau =2$ seems to be optimal for defining $\B(m,A)$.\\

{\bf Proof of Theorem~\ref{theo1.2}.} Part (ii) of the claim immediately follows from \cite[Lemma~6.1]{win_bu}. In fact,
  given $m>0$, $p \in (1, \frac{2n}{n+2})$ as well as radial and positive functions $u \in C^0(\bar\Omega)$ and $v\in 
  W^{1,\infty}(\Omega)$ with $\io u =m$, sequences $(u_k)_{k\in\N} \subset C^0(\bar\Omega)$
  and $(v_k)_{k\in\N} \subset W^{1,\infty}(\Omega)$ of radially symmetric positive functions with $\io u_k=m$ for all
  $k\in\N$ are constructed, which satisfy
  \be{1.2.2}
	u_k\to u \quad \mbox{in } L^p(\Omega),
	\quad
	v_k \to v \quad \mbox{in } W^{1,2}(\Omega), \quad\mbox{and}\quad  \io u_k v_k \to \infty
	\qquad \mbox{as } k\to\infty.
  \ee
  Combining this with \eqref{G1} and our additional condition $p > 2-\alpha$, we find some $C>0$ such that 
  $$\frac{1}{2} \io |\nabla v_k|^2 + \frac{1}{2} \io v_k^2 + \io G(u_k) \le C \qquad\mbox{for all } k \in \mathbb{N}.$$ 
  Thus, \eqref{1.2.2} implies $\F(u_k,v_k) \to -\infty$ as $k \to \infty$ which proves part (ii) of the claim. \\

  In view of part (ii), it is sufficient to prove part(i) of the claim in the case $\alpha \in (\frac{2}{n},1)$. To this end we notice that,
  given $m>0$ and
  \begin{equation}\label{1.2.3}
   \gamma_2 \in ((1-\alpha)n, n-2),     
  \end{equation}  
  by \cite[Lemma~4.1]{win_mmas} there exists $\eta_0 >0$ such that for any $\eta \in (0,\eta_0)$ there are radial and positive 
  functions $u_\eta, v_\eta \in C^\infty (\bar{\Omega})$ with $\io u_\eta =m$ satisfying 
  $$\io |\nabla v_\eta|^2 \le c_1 \eta^{-(-n +2\gamma_2+2)}, \;
    \io v_\eta^2 \le c_1 \eta^{-(-n +2\gamma_2)}, \;
    \io G(u_\eta) \le c_1 \eta^{-(1-\alpha)n}, \;
    \io u_\eta v_\eta \ge c_2 \eta^{-\gamma_2}$$
  for all $\eta \in (0,\eta_0)$ with positive constants $c_1$ and $c_2$. Hence, \eqref{1.2.3} implies that there are 
  $c_3, c_4 >0$ and $\eta_1 \in (0,\eta_0)$ such that   
  $$\| v_\eta \|_{W^{1,2} (\Omega)} \le A_\eta :=c_3 \eta^{-(\gamma_2+1-\frac{n}{2})} \qquad\mbox{and}\qquad 
    \F(u_\eta,v_\eta) \le -c_4 \eta^{-\gamma_2} \qquad\mbox{for all } \eta \in (0,\eta_1)$$
  are fulfilled. Since $\gamma_2 <n-2$ implies $\gamma_2 > 2 (\gamma_2+1-\frac{n}{2})$, we conclude that there exist $\eta_2 \in
  (0,\eta_1)$ and $c_5 >0$ such that
  $$\F(u_\eta,v_\eta) \le - K(m) \left( 1+ A_\eta^2 \right) \qquad\mbox{for all } \eta \in (0,\eta_2).$$
  Hence, $(u_\eta, v_\eta) \in \B (m,A_\eta)$ for $\eta$ small enough.  
\qed

\mysection{Unbounded global-in-time solutions}\label{section5}

The last section is devoted to the proof of Theorem \ref{theo1.5}. To this end we provide the following lemma.

\begin{lem}\label{lemat}
Let $\Omega\subset \R^n$ with some $n\ge 2$ and $R>0$. Moreover assume that \eqref{balance} and \eqref{decay}
are satisfied.
Then there exists $p>n$ such that for any solution $(u,v)$ to \eqref{0} and any $T \in (0, \infty)$ with $T \le \tm$ there is $C>0$ such that $u$ admits the estimate
\begin{equation}\label{wzor}
\|u(\cdot,t)\|_{L^p(\Omega)}\leq C, \qquad t \in \left( \frac{T}{2},T \right).
\end{equation}
\end{lem}  
Before proving the above lemma let us show how to infer Theorem \ref{theo1.5} from it. 

{\bf Proof of Theorem \ref{theo1.5}.} We fix $T \in (0, \infty)$ with $T \le \tm$ and first use the second equation in \eqref{0}. By a standard regularity result in the theory of parabolic equations, see \cite[Lemma 4.1]{horstmann_winkler} for example, (\ref{wzor}) yields a uniform estimate of the $L^\infty$-norm of $\nabla v$ on $(\frac{T}{2},T)$. Then by (\ref{balance}) and \cite[Theorem 2.2]{ja} we arrive at the uniform estimate of $\|u\|_{L^\infty(\Omega)}$ on $(\frac{T}{2},T)$. Hence,  in view of Lemma \ref{lem2.1}, we have shown the existence of a global-in-time solution to \eqref{0} whatever initial data we start with. On the other hand choosing $\Omega=B_R$ and radially symmetric initial data, since (\ref{H1}) and (\ref{G1}) are satisfied, we conclude with the use of \cite[Theorem~5.1]{win_mmas} that the  solutions we arrived at are unbounded.
\qed    

Next we complete this section by proving Lemma \ref{lemat}.

\noindent 
{\bf Proof of Lemma \ref{lemat}.} Multiplying the first equation of \eqref{0} by $u^{p-1}$, $p \in (n, \gamma_1]$, and the second one by
$\Delta v$, we arrive at
\begin{equation}
\label{eins}
\frac{1}{p}\frac{d}{dt}\int_\Omega u^pdx + (p-1)\int_\Omega \phi(u)\left|\nabla u\right|^2u^{p-2}dx = (p-1)\int_\Omega u^{p-1}\beta(u)\nabla v\nabla u \; dx,
\end{equation}
and 
\begin{equation}
\label{zwei}
\frac{1}{2}\frac{d}{dt}\int_\Omega |\nabla v|^2dx + \frac{1}{2}\int_\Omega\left|\Delta v\right|^{2}dx+ \int_\Omega|\nabla v|^2dx\leq \frac{1}{2}\int_\Omega u^2dx.
\end{equation}
Since
\[
u^{p-1}\beta(u)=u^{\frac{p-2}{2}}u^{\frac{p}{2}}\sqrt{\beta(u)}\sqrt{\beta(u)},
\]
in view of (\ref{balance}) we infer from (\ref{eins}) that
\begin{equation}\label{drei}
\frac{1}{p}\frac{d}{dt}\int_\Omega u^pdx + \frac{p-1}{2}\int_\Omega \phi(u)\left|\nabla u\right|^2u^{p-2}dx\leq
C\int_\Omega u^p\beta(u)|\nabla v|^2dx.
\end{equation}
Next adding (\ref{drei}) and (\ref{zwei}) and applying (\ref{decay}) we arrive at
\begin{equation}\label{vier}
\frac{d}{dt}\left(\int_\Omega u^pdx+\int_\Omega|\nabla v|^2dx\right)\leq C\left(\int_\Omega u^pdx\right)^{\frac{2}{p}}+C\int_\Omega|\nabla v|^2\leq C\left(\int_\Omega u^pdx+\int_\Omega|\nabla v|^2dx +1\right),
\end{equation}
which in turn, by Gr\"{o}nwall's lemma yields the claimed estimate of $\|u\|_{L^p(\Omega)}$. 

\qed

\end{document}